\documentclass[a4paper,12pt,final]{amsart}
\usepackage{times,a4wide,mathrsfs, nameref,showkeys,amssymb,amsmath,amsthm,enumerate,xypic,tikzsymbols,dsfont}

\newcommand{\C}{\mathbb{C}}

\newcommand{\QQ}{\mathbb{Q}}
\newcommand{\NN}{\mathbb{N}}
\newcommand{\PP}{\mathbb{P}}

\newcommand{\Ss}{\mathcal S}

\newcommand{\XX}{\mathcal X}
\newcommand{\YY}{\mathcal Y}

\newcommand{\Zz}{\mathcal Z}

\newcommand{\MM}{\mathcal M}

\newcommand{\FF}{\mathcal F}

\newcommand{\pic}{\hbox{Pic}}

\newcommand{\wt}{\widetilde}
\newcommand{\rom}{\romannumeral}

\newcommand{\one}{\mathds{1}}

\DeclareMathOperator{\aut}{Aut}
\DeclareMathOperator{\bir}{Bir}
\DeclareMathOperator{\ide}{id}

\DeclareMathOperator{\ima}{Im}

\DeclareMathOperator{\sym}{Sym}
\DeclareMathOperator{\Gr}{Gr}

\DeclareMathOperator{\id}{id}

\newtheorem{theorem}{Theorem}[section]

\newtheorem{lemma}[theorem]{Lemma}

\newtheorem{corollary}[theorem]{Corollary}
\newtheorem{proposition}[theorem]{Proposition}
\newtheorem{conjecture}[theorem]{Conjecture}
\newtheorem{remark}[theorem]{Remark}
\newtheorem{definition}[theorem]{Definition}
\newtheorem{convention}{Conventions}

\newtheorem{notation}[theorem]{Notation}

\newtheorem{nonumbering}{Theorem}

\newtheorem{nonumberingp}{Proposition}

\newtheorem{nonumberingt}{Acknowledgements}

\newif\ifHideFoot
\HideFoottrue  
\HideFootfalse  

\ifHideFoot

\newcommand{\Robert}[1]{}
\newcommand{\Michele}[1]{}

\else 

\newcommand{\marg}[1]{\normalsize{{
			\color{red}\footnote{{\color{blue}#1}}}{\marginpar[\vskip
			-.25cm{\color{red}\hfill$\Rightarrow$\tiny\thefootnote}]{\vskip
				-.2cm{\color{red}$\Leftarrow$\tiny\thefootnote}}}}}
\newcommand{\Robert}[1]{\marg{(Robert) #1}}
\newcommand{\Michele}[1]{\marg{(Michele) #1}}

\fi

\begin{document}

\author[Michele Bolognesi]{Michele Bolognesi}
\address{Institut Montpellierain Alexander Grothendieck \\ %
Universit\'e de Montpellier \\ %
CNRS \\ %
Case Courrier 051 - Place Eug\`ene Bataillon \\ %
34095 Montpellier Cedex 5 \\ %
France}
\email{michele.bolognesi@umontpellier.fr}

\author[Robert Laterveer]
{Robert Laterveer}

\address{Institut de Recherche Math\'ematique Avanc\'ee,
CNRS -- Universit\'e 
de Strasbourg,\
7 Rue Ren\'e Des\-car\-tes, 67084 Strasbourg CEDEX,
FRANCE.}
\email{robert.laterveer@math.unistra.fr}

\title{Double EPW sextics and the Voisin filtration on zero-cycles}

\begin{abstract} Let $X$ be a double EPW sextic, and $\iota$ its anti-symplectic involution. We relate the $\iota$-anti-invariant part of the Chow group of zero-cycles of $X$
with Voisin's rational orbit filtration. For a general double EPW sextic $X$, we also relate the anti-invariant part of the Chow motive of $X$ with the motive of a Gushel--Mukai fourfold.
As an application of the first result, we obtain a similar result for certain Fano varieties of lines in cubics with infinite-order birational automorphisms.
\end{abstract}

\thanks{\textit{2020 Mathematics Subject Classification:} 14C15, 14C25, 14J42}
\keywords{Chow groups, motive, hyper-K\"ahler varieties, Bloch conjecture, Beauville ``splitting property'' conjecture, Voisin's rational orbit filtration}
\thanks{MB and RL are supported by ANR grant ANR-20-CE40-0023.}

\maketitle

\section{Introduction}

Let $X$ be a smooth projective variety of dimension $n$. Bloch and Beilinson have conjectured that there should exist a decreasing filtration $F_{BB}^\ast$ on the Chow groups $A^\ast(X)$, satisfying a number of axioms \cite{Jan}. In the special case that $X$ is a hyper-Kähler variety, one can (conjecturally) say more. In this case, Beauville \cite{Beau3} has conjectured that the Chow ring should split into a bigraded ring

\begin{equation}\label{abx} A^{\ast}(X)= A^{\ast}_{(\ast)}(X)\ ,\end{equation}
with $A^{j}_{(i)}(X)= F_{BB}^{i} A^j(X)/F_{BB}^{i+1} A^j(X)$. 

To try and make sense of this, Voisin \cite{Vcoiso} has defined a filtration $S_\ast$ on the Chow group of zero-cycles of a hyper-K\"ahler variety $X$ of dimension $n=2m$. In particular, she expects this filtration to be opposite to the Bloch-Beilinson filtration (for more on what it means for a filtration on zero-cycles to be opposite to the conjectural filtration $F^\ast_{BB}$ cf. \cite{Vcoiso} and \cite{vialbirat}). The Voisin filtration on zero-cycles is defined as follows:

$$S_jA^{n}(X):= \langle x\in X| \dim O_x \geq m- j \rangle\ \ \subset\ A^n(X)\ ,$$
where $O_x\subset X$ is the locus of points rationally equivalent to $x$. 
Assuming (as in \cite{Vcoiso}) that $S_\ast$ and $F^{2\ast}$ are opposite filtrations on $A^n(X)$, there should be canonical isomorphisms

\begin{equation}\label{split} S_{i} / S_{i-1} \xleftarrow{\cong} S_{i}\cap F_{BB}^{2m-2i} \xrightarrow{\cong} F_{BB}^{2m-2i}/F_{BB}^{2m-2i+2}=:  A^n_{(n-2i)}(X)\ .
\end{equation}


Let us now assume that $X$ has a birational automorphism $\phi$. The axioms of the Bloch--Beilinson filtration then imply that $\phi$ respects the bigrading \eqref{abx}, and that

$$\phi^{\ast}:A^n_{(2i)}(X) \to A^n_{(2i)}(X)$$
should be equal to $(-1)^i \id$ if $\phi$ is anti-symplectic and equal to $\id$ if $\phi$ is symplectic. 
Combined with the conjectural equality \eqref{split}, this yields the following conjecture:

\begin{conjecture}\label{conj} Let $X$ be a hyper-K\"ahler variety of dimension $n$, and $\phi\in\bir(X)$ a birational automorphism.

\noindent
(\rom1) If $\phi$ is symplectic then 
  \[ \phi^\ast=\id\colon\ A^n(X)\ \to\ A^n(X)\ . \]
  
 \noindent
 (\rom2) If $\phi$ is anti-symplectic then ($\phi$ respects $S_\ast$ \footnote{This is clearly the case when $\phi$ is a finite-order automorphism.} and)
   \[ \phi^\ast = (-1)^i \id\colon\ \   \Gr_{i}^S A^n(X) \ \to\  \Gr_{i}^S A^n(X)\ .\]
 \end{conjecture}
 
 This conjecture has been established in specific instances: the conjecture is true for $n=2$ and $\phi\in\aut(X)$ of finite order \cite{Huy}. Further
 cases of Conjecture \ref{conj}(\rom1) have been established in \cite{Fu}, \cite{vial}, \cite{LYZ}, \cite{72}, while cases of Conjecture \ref{conj}(\rom2) are proven in \cite{25}, \cite{26}, \cite{28}, \cite{27}, \cite{LV}, \cite{61}, \cite{65}, \cite{LYZ}, \cite{72}. 
 


In this paper, we consider Conjecture \ref{conj} for double EPW sextics. These varieties form a locally complete family of hyper-K\"ahler varieties of $K3^{[2]}$-type, and by construction they come with an anti-symplectic involution. Our main result establishes Conjecture \ref{conj} in new cases:

\begin{nonumbering}[= Theorem \ref{mainequality}]

Let $X$ be a double EPW sextic and  $\imath: X \to X$ the natural anti-symplectic involution, and let $S_\ast A^4(X)$ denote Voisin's filtration on zero-cycles. Then we have an identification 

\begin{equation}\label{mainequation}
A^4(X)^-=S_1 A^4(X)\cap A^4_{hom}(X)\ \end{equation}
(where $A^4(X)^-$ denotes the $(-1)$-eigenspace with respect to the action of $\iota$).

\end{nonumbering}

The existence of this identity had been conjectured by Li--Yu--Zhang in \cite[Section 3]{LYZ}. Theorem \ref{mainequality} generalizes \cite[Theorem 2]{LV}, where this was proven when $X$ is birational to $S^{[2]}$, for $S$ a K3 surface. Theorem \ref{mainequality} implies in particular that $S_1 A^4(X)$ is motivic, i.e. it is cut out by a projector (Corollary \ref{motivicS1}).

More precisely, the equality \eqref{mainequation} comes from an eigen-space decomposition

\begin{equation}\label{eigen} A^4(X)= \mathbb Q [o_X] \oplus A^4(X)_{}^- \oplus A^4_{hom}(X)^+_{}\ ,\end{equation}
where $o_X$ is any point of the fixed locus surface of $\imath$ and $\imath_*$ operates as $\ide$ (resp. $-\ide$) on
$A^4_{hom}(X)^+_{}$  (resp. on $A^4(X)^-_{}$). 

Our result confirms the expectation that the splitting \eqref{eigen} is the ``right'' motivic splitting for zero-cycles on double EPW sextics, i.e. that it should be the restriction to $A^4(X)$ of the conjectural splitting \eqref{abx} put forth by Beauville. Indeed, as a consequence of our main result, we obtain identifications

    \begin{align*}
  &S_1 A^4(X) =  \mathbb Q [o_X] \oplus A^4(X)_{}^- = A^1(X)\cdot A^3(X) = A^2(X)\cdot h^2 \ ,\\
  & S_0 A^4(X)= \QQ[o_X]=\QQ[h^4]=\QQ[c_4(T_X)]= S_1 A^4(X)\cap A^4(X)^+\\
  \end{align*}
(where $h\in A^1(X)$ denotes the polarization), cf. Remark \ref{stuff} and Section 4 below.
The different characterizations of $S_1 A^4(X)$ given here are similar to known results for Fano varieties of lines on cubic fourfolds (cf. \cite{SV}). The second line gives a new characterization for what should be the ``distinguished 0-cycle'' on $X$, as any degree 1 zero-cycle that is $\iota$-invariant and supported on a uniruled divisor.


When $X$ is a general double EPW sextic, there is yet another characterization of $S_1 A^4(X)$: a result of Zhang \cite{Zh} implies that there is an isomorphism
  \[  \Gamma_\ast \colon\ A^3(Y) \ \xrightarrow{\cong}\ S_1 A^4(X)\ ,\]
  where $Y$ is a Gushel--Mukai fourfold and $\Gamma$ is some natural correspondence defined in terms of conics contained in $Y$ (cf. Proposition \ref{0zh} below); again, this is in line with known results about the Chow ring of Fano varieties of lines on cubic fourfolds \cite{SV}.


As a further result of independent interest, we establish the following relation of motives:

\begin{nonumberingp}[=Proposition \ref{isomotives}] Let $X$ be a general double EPW sextic, so that $X$ is associated to a Gushel--Mukai fourfold $Y$ via the construction of Iliev--Manivel (cf. Theorem \ref{im} below). Let us write $X\to Z$ for the double cover, where $Z\subset\PP^5$ is an EPW sextic. There is an isomorphism of Chow motives
  \[  h(X) \cong h(Z) \oplus \bigoplus_{i=-1}^1 t(Y)(i)\oplus \bigoplus_{j=-3}^{-1} \one(j)^{\oplus r_j}\ \ \ \hbox{in} \ \MM_{\rm rat}\ ,\]
 where $t(Y)$ denotes the transcendental motive of $Y$. 
\end{nonumberingp}

This builds on, and improves, the afore-mentioned result of Zhang about zero-cycles on $X$.

\medskip

Using recent results of Brooke--Frei--Marquand--Qin \cite{BFMQ}, our Theorem \ref{mainequality} also allows to prove Conjecture \ref{conj} for the Fano variety of lines of a very general cubic fourfold in the Hassett divisor $\mathcal{C}_{12}$:

\begin{nonumberingp}[= Proposition \ref{conjC12}]

Let $Y$ be a very general cubic fourfold in the divisor $\mathcal{C}_{12}$, let $X=F(Y)$ be its Fano variety of lines and let $\phi$ be any birational automorphism of $X$. Then Conjecture \ref{conj} holds true for $(X,\phi)$. 

\end{nonumberingp} 

One reason this case is particularly interesting is that Hassett--Tschinkel \cite{HT} have shown that varieties $X$ as in Proposition \ref{conjC12} admit (symplectic and non-symplectic) birational automorphisms of infinite order.

\begin{convention} In this article, the word {\sl variety\/} will mean a reduced irreducible scheme of finite type over $\C$. A {\sl subvariety\/} will refer to a (possibly reducible) reduced subscheme which is equidimensional. 

{\bf All Chow groups will be with rational coefficients}: we denote as customary by $A_j(Y)$ the Chow group of $j$-dimensional cycles on $Y$ with $\QQ$-coefficients; for $Y$ a smooth variety of dimension $n$ the notations $A_j(Y)$ and $A^{n-j}(Y)$ are used interchangeably. 
The notation $A^j_{hom}(Y)$ will be used to denote the subgroup of homologically trivial cycles.
For a morphism $f\colon X\to Y$, we write $\Gamma_f\in A_\ast(X\times Y)$ for the graph of $f$.

The contravariant category of Chow motives (i.e., pure motives with respect to rational equivalence as described in \cite{Sc}, \cite{MNP}) will be denoted 
$\MM_{\rm rat}$. We will write $H^j(X)$ to indicate singular cohomology $H^j(X,\QQ)$.
\end{convention}

 \section{Preliminaries}

 \subsection{Double EPW sextics: the original construction}
 
Double EPW sextics are double covers of certain sextic hypersurfaces:
	
	\begin{definition}[Eisenbud--Popescu--Walter \cite{EPW}] 
		Let $A\subset \wedge^3 \C^6$ be a subspace which
		is Lagrangian with respect to the symplectic form on $\wedge^3 \C^6$ given by
		the wedge product. Define the degeneracy locus		
		\[ Z_A:= \Bigl\{  [v]\in \PP(\C^6)\ \vert\ \dim \bigl( A\cap ( v\wedge
		\wedge^2 \C^6)\bigr) \ge 1\Bigr\}\ \ \subset \PP(\C^6)\  .\]
		An {\em EPW sextic\/} is a  $Z_A$ for some Lagrangian subset $A\subset \wedge^3 \C^6$.
	\end{definition}

	\begin{theorem}[O'Grady]
		Let $Z=Z_A$ be an EPW sextic that is a hypersurface and such that the singular locus $S:=\hbox{Sing}(Z)$ is a
		smooth irreducible surface. Let 
		  \[ f\colon\ X\ \to\  Z\] 
		  be the double cover branched over
		$S$. Then $X$
		is a smooth hyper-K\"ahler fourfold of K3$^{[2]}$-type (a so-called {\em double
			EPW sextic}), 
		and the class 
		  \[h:=f^*c_1(O_Z(1))\ \ \in \ A^1(X)\]
		  defines a polarization of square
		$2$ for the Beauville--Bogomolov form.
		Double EPW sextics form a $20$-dimensional locally complete family.
			\end{theorem}
	
	\begin{proof} 
		This is \cite[Theorem 1.1(2)]{OG2}. We observe that the hypothesis on $Z$ and
		$\hbox{Sing}(Z)$ are satisfied by a generic EPW sextic (indeed, it
		suffices that the Lagrangian subspace $A$ be in $\hbox{LG}(\wedge^3 V)^0$, which
		is a certain open dense subset of $\hbox{LG}(\wedge^3 V)$ defined in
		\cite[Section 2]{OG2}). Letting $A$ vary in $\hbox{LG}(\wedge^3 V)^0$, we
		obtain a locally complete family with $20$ moduli (as observed in
		\cite[Introduction]{OG2}).		
	\end{proof}

 \subsection{The Chow ring of double EPW sextics}
 
 Thanks to the work of the second named author (and building on earlier work by Ferretti \cite{Fe0}, \cite{Fe} and by Laterveer--Vial \cite{LV}), the Beauville--Voisin conjecture is now settled for double EPW sextics:
 
\begin{theorem}[\cite{Lbv}]
Let $X$ be a smooth double EPW sextic. The $\QQ$-subalgebra

$$\langle A^1(X),c_j(X)\rangle \subset A^\ast (X)$$
injects into cohomology, under the cycle class map.

In addition, let $A^2(X)^+\subset A^2(X)$ denote the subgroup of cycles invariant under the covering involution $\iota$. The cycle class map induces injections

$$\langle A^1(X),c_j(X),A^2(X)^+ \rangle \cap A^i(X) \hookrightarrow H^{2i}(X,\QQ),\ for\ i\geq 3.$$
\end{theorem}

\begin{proof} This is \cite[Theorem 3.1]{Lbv}.
\end{proof}
 
 
   
   
   
  
  
  \subsection{Double EPW sextics: the Iliev--Manivel construction}
  
The general double EPW sextic can alternatively be constructed in terms of conics on a Gushel--Mukai fourfold.
  
  \begin{definition} An {\em ordinary Gushel--Mukai fourfold\/} is a smooth dimensionally transverse complete intersection
    \[ Y:= \Gr(2,5)\cap H\cap Q\ \ \ \subset\ \PP^9\ ,\]
    where $\Gr(2,5)$ is the Grassmannian of 2-dimensional subspaces of a fixed 5-dimensional complex vector space and $H, Q$ are a hyperplane resp. a quadric with respect to the Pl\"ucker embedding.
    \end{definition}
    
    Gushel--Mukai varieties have been studied in depth by Debarre and Kuznetzov \cite{D}, \cite{DK}, \cite{DK1}, \cite{DK3}, \cite{DK2}, \cite{DK4}.

    \begin{notation} Let us call $F(G)$ the Hilbert scheme parametrizing conics contained in $\Gr(2,5)$. Given an ordinary Gushel--Mukai fourfold $Y$, let $F=F(Y)$ denote the Hilbert scheme of conics contained in $Y$. We write $P\subset F\times Y$ for the universal conic, with projections $p\colon P \to F$ and $q\colon P\to Y$.
    
Conics in $\Gr(2,5)$ can be of type $\tau$ or of type $\rho$ or of type $\sigma$ (cf. \cite[Section 3.1]{IM}); we write
    \[ F(Y)= F_\tau(Y)\cup F_\rho(Y)\cup F_\sigma(Y) \]
    for the partition of $F(Y)$ according to the type of conic.
      \end{notation}
      
      \begin{theorem}(Iliev--Manivel \cite{IM})\label{im} Let $Y$ be a general Gushel--Mukai fourfold. The variety $F=F(Y)$ is a 5 dimensional smooth projective variety, and the subvarieties $F_\rho(Y), F_\sigma(Y)\subset F$ are of codimension 1 resp. 2. There exists a natural morphism
        \[ \pi\colon F\ \to\ X\ ,\]
        where $X$ is a double EPW sextic. The morphism $\pi$ restricted to $F_\tau(Y)$ is a $\PP^1$-fibration, and $\pi$ contracts $F_\rho(Y)$ and $F_\sigma(Y)$ to points $x_\rho$ resp. $x_\sigma$ on $X$.
        
        This construction defines a dominant rational map 
        \[\nu\colon\ \  \MM_{GM4}\ \dashrightarrow\ \MM_{dEPW} \]
        from the moduli space of Gushel--Mukai fourfolds to the moduli space of double EPW sextics.
        \end{theorem}
        
        \begin{proof} The smoothness of $F$ is contained in \cite[Theorem 3.2]{IM}; the relation between $F$ and a double EPW sextic is \cite[Proposition 4.18]{IM}, and the fact that a general double EPW sextic is attained in this way comes from \cite[Corollary 4.17]{IM}.
         \end{proof}
         
   \begin{remark} Theorem \ref{im} has recently been reproven (and made more precise) in \cite[Theorem 7.12]{DK4}.
    \end{remark}

\subsection{Double EPW sextics: the modular construction}    

The relation between the very general double EPW sextic and a Gushel--Mukai fourfold has been recently described in terms of moduli spaces:

\begin{theorem}(Perry--Pertusi--Zhao \cite{PPZ}, Guo--Liu--Zhang \cite{GLZ})\label{ppz} Let $Y$ be a very general Gushel--Mukai fourfold. There exists a Mukai vector $v$ and stability conditions $\sigma$ on the Kuznetsov component $Ku(Y)$, such that
there is an isomorphism
  \[ X\cong M_{v,\sigma}(Ku(Y)) \]
  between a double EPW sextic $X$ and the moduli space $M_{v,\sigma}(Ku(Y))$. The very general double EPW sextic is attained in this way.
  
  (That is, there exists $\MM^\circ_{GM4}\subset\MM_{GM4}$, intersection of a countable infinity of dense open subsets, such that over $\MM^\circ_{GM4}$ the double EPW sextic has a modular interpretation.)
  
  \end{theorem}    

\begin{proof} Thanks to
\cite[Theorem 1.7]{PPZ}, for a general $Y$, the moduli space  $M_{v,\sigma}(Ku(Y))$ exists and is a smooth projective hyper-K\"ahler variety. For the very general $Y$, the moduli space coincides with a double EPW sextic as shown in \cite[Proposition 5.17]{PPZ}. 
   
For later use, we observe that the construction of \cite{PPZ} is actually done in families over the moduli space of Gushel--Mukai fourfolds.
\end{proof}

\begin{remark} 

While we will not need this, we observe that in \cite[page 30]{F+} it is asserted that in Theorem \ref{ppz} the words ``very general'' can be replaced by ``general''.
  \end{remark}

As a consequence of the modular construction, we can relate double EPW sextics and Gushel--Mukai fourfolds on the level of Chow motives:

\begin{theorem}\label{bul} Let $Y$ be a very general Gushel--Mukai fourfold, and let $X$ be the double EPW sextic associated to $Y$ via Theorem \ref{ppz}. There is an inclusion of Chow motives
  \[  h(X)\ \hookrightarrow\ \bigoplus h(Y^2)(\ast)\ \ \ \hbox{in}\ \MM_{\rm rat}\ .\]
 \end{theorem}
 
 \begin{proof} This comes essentially from B\"ulles' result \cite{Bul}, with the improved bound on the exponent obtained in \cite[Theorem 1.1]{FLV3}. For later use, we remark that the argument proving Theorem \ref{bul} can be done in families over a base; this is explained in \cite{FLV3}.
 \end{proof}

   \subsection{Generically defined cycles and the Franchettina property}  
   \label{ss:gen}
   
   \begin{definition}\label{frank} Let $\YY\to S$ be a smooth projective morphism, where $\YY, S$ are smooth quasi-projective varieties. We will say that $\YY\to S$ has the {\em Franchetta property in codimension $j$\/} if the following holds: for every $\Gamma\in A^j(\YY)$ such that the restriction $\Gamma\vert_{Y_s}$ is homologically trivial for all $s\in S$, the restriction $\Gamma\vert_s$ is zero in $A^j(Y_s)$ for all $s\in S$.
 
 We say that $\YY\to S$ has the {\em Franchetta property\/} if $\YY\to S$ has the Franchetta property in codimension $j$ for all $j$.
 \end{definition}
 
 This property is studied in \cite{BL}, \cite{FLV}, \cite{FLV3}, \cite{FLV2}. The following useful ``spread lemma'' implies that in Definition \ref{frank}, it actually suffices to consider the very general fiber:
 
 \begin{lemma}\label{spread} Let $\YY\to S$ be a smooth projective morphism, and $\Gamma\in A^j(\YY)$. Assume that the restriction $\Gamma\vert_{Y_s}$ is 0 in $A^j(Y_s)$ for the very general $s\in S$. Then the restriction $\Gamma\vert_{Y_s}$ is 0 in $A^j(Y_s)$ for all $s\in S$. 
 \end{lemma}
 
 \begin{proof} This is \cite[Lemma 3.1]{Vo}. A proof can be found in \cite[Proposition 2.4]{V3}.
  \end{proof}

 We can also consider a weaker variant of the Franchetta property, by replacing Chow groups $A^\ast()$ by groups $B^\ast()$ of algebraic cycles modulo algebraic equivalence:
 
  \begin{definition}\label{frank2} Let $\YY\to S$ be a smooth projective morphism, where $\YY, S$ are smooth quasi-projective varieties. We will say that $\YY\to S$ has the {\em Franchettina property in codimension $j$\/} if the following holds: for every $\Gamma\in B^j(\YY)$ such that the restriction $\Gamma\vert_{Y_s}$ is homologically trivial for all $s\in S$, the restriction $\Gamma\vert_s$ is zero in $B^j(Y_s)$ for all $s\in S$.
 
 We will say that $\YY\to S$ has the {\em Franchettina property\/} if $\YY\to S$ has the Franchettina property in codimension $j$ for all $j$.
 \end{definition}
 
 The Franchettina property is studied in \cite{BoLa}.

 \begin{notation}\label{def:gd} Given a family $\YY\to S$ as above, with $Y:=Y_s$ a fiber, we write
   \[ \begin{split} GDA^j_S(Y):=\ima\Bigl( 
  A^j(\YY)\to A^j(Y)\Bigr)\ ,\\
      GDB^j_S(Y):=\ima\Bigl( 
  B^j(\YY)\to B^j(Y)\Bigr)\ \\   
  \end{split}\]
   for the subgroups of {\em generically defined cycles}. 
Whenever it is clear to which family we are referring, the index $S$ will often be suppressed from the notation.
  \end{notation}
  
  With this notation, the Franchetta (and Franchettina) property amounts to saying that $GDA^\ast_S(Y)$ (resp. $GDB^\ast_S(Y)$) injects into cohomology, via the cycle class map. 
  
  Let us recall the following result, which we will need later on:
   
   \begin{proposition}\label{franche} Let $\XX^\prime\to\MM^1_{GM4}$ be the universal family of hyper-K\"ahler varieties as constructed in \cite{PPZ} (cf. Theorem \ref{ppz}), where $\MM^1_{GM4}$ is an open subset of $\MM_{GM4}$. This family has the Franchettina property, i.e.
   the cycle class map induces injections
       \[ GDB^\ast_{\MM^1_{GM4}}(X)\ \to\ H^\ast(X,\QQ) \ ,\]
       for the very general fiber $X$.   
   \end{proposition}
   
   \begin{proof} 
This is Proposition 2.15 from \cite{Lbv} (see also Corollary 8.2 in \cite{BoLa}, where this is proven for $ GDB^\ast_{\MM_{dEPW}}(X)$, where $ \MM_{dEPW}$ denotes the moduli space of double EPW sextics). This follows from Theorem \ref{bul}, in combination with the fact that one can readily prove the Franchettina property for the fiber square $\YY^{2/\MM_{GM4}}\to \MM_{GM4}$, where $\YY\to \MM_{GM4}$ is the universal GM fourfold.
\end{proof}

In this paper, the Franchettina result will be used in the following guise:  
  
  \begin{corollary}[\cite{Lbv}]\label{Fr2} Let $\XX\to \MM_{GM4}^2$ be the universal family of hyper-K\"ahler varieties as constructed in Theorem \ref{im}, where $\MM^2_{GM4}$ is a dense open subset of $\MM_{GM4}$. 
 For the very general fiber $X$ of $\XX\to \MM_{GM4}^2$, we have
    \[ \begin{split}  
    GDB^2_{\MM^2_{GM4}}(X) &=  \QQ[h^2]\oplus \QQ[c_2]\ ,\\
     GDB^3_{\MM^2_{GM4}}(X) &=  \QQ[h^3]\ ,\\  
     \end{split}  \]
where $h$ is the linear section and $c_2$ the second Chern class of $X$.
 \end{corollary}
 
 \begin{proof} Over the common open $\MM^{12}_{GM4}:=\MM^1_{GM4}\cap \MM^2_{GM4}$, there is a morphism from $\XX$ to $\XX^\prime$ which induces an isomorphism on the very general fiber by Theorem \ref{ppz} (as explained in \cite{PPZ} and \cite{GLZ}, the morphism from $\XX$ to $\XX^\prime$ is induced by the morphism sending a conic $c$ in the GM fourfold $Y$ to $pr_Y(I_c)$, where $I_c$ is the ideal sheaf of $c$ and $pr_Y$ is a projection functor; this is defined relatively over the open $\MM_{GM4}^{12}$). Thus, generically defined cycles for the family $\XX$ correspond to generically defined cycles for the family $\XX^\prime$, and so the Franchettina property for the family $\XX$ follows from Proposition \ref{franche}.
 A direct identification of the Hodge classes in codimension 2 (cf. \cite[Corollary 2.20]{Lbv}) then yields the Corollary.
 \end{proof}
  
We observe that Corollary \ref{Fr2} corrects the statement of \cite[Proposition 2.15]{Lbv}, where the restriction of the family to a countable union of dense open subsets is considered.
  
 
 \subsection{Two families} There are 2 (fiberwise isomorphic, but a priori different) families of double EPW sextics over $\MM^2_{GM4}$: there is the family $\XX\to\MM^2_{GM4}$ constructed by Iliev--Manivel (cf. Theorem \ref{im}), and there is the base changed family $\XX^\nu:= \XX\times_{\MM_{dEPW}} \MM^2_{GM4}$ (where $\XX\to\MM_{dEPW}$ is the universal family over the moduli space of double EPW sextics). For later use, we establish a lemma comparing these two families:
 
 \begin{lemma}\label{2fam} Let $X$ and $X^\nu$ be fibers of the family $\XX\to\MM^2_{GM4}$ resp. $\XX^\nu\to\MM^2_{GM4}$ over the same point. Let $f\colon X\to X^\nu$ be the isomorphism between the two fibers.
 
 \noindent
 (\rom1) The isomorphism $f_\ast\colon A^\ast(X)\xrightarrow{\cong} A^\ast(X^\nu)$ induces an isomorphism 
     \[  f_\ast\colon GDA^\ast_{\MM^2_{GM4}}(X)\ \xrightarrow{\cong}\  GDA^\ast_{\MM^2_{GM4}}(X^\nu) \ .\]
     
 \noindent
 (\rom2)      The isomorphism $A^\ast(X\times X)\xrightarrow{\cong} A^\ast(X^\nu\times X^\nu)$ induces an isomorphism 
     \[  (f,f)_\ast\colon GDA^\ast_{\MM^2_{GM4}}(X\times X)\ \xrightarrow{\cong}\  GDA^\ast_{\MM^2_{GM4}}(X^\nu\times X^\nu) \ .\]
 \end{lemma}
  
  \begin{proof} Let us prove (\rom1); the argument for (\rom2) is only notationally more complicated.
  
  We need to understand the construction of the isomorphism $f\colon X\to X^\nu$ in \cite{IM}. This construction involves a novel interpretation of EPW sextics: as in \cite[Section 2]{IM}, let $I^\vee\cong\PP^5$ denote the projective dual of the linear system $I$ of quadrics containing a given GM fourfold $Y$. As shown in loc. cit., there is a hypersurface $Z\subset I^\vee$ (defined in terms of quadrics for which the restriction to a certain subspace is singular; $Z$ is denoted $Y^\vee_Z$ in \cite{IM}) that is isomorphic to an EPW sextic.
 This construction can be done relatively over $\MM^2_{GM4}$; the relative version of $I^\vee$ is isomorphic to $\PP^5\times\MM^2_{GM4}$ and this isomorphism restricts to give an isomorphism of families of EPW sextics over $\MM^2_{GM4}$
    \[ g\colon\ \ \  \Zz\ \xrightarrow{\cong}\ \Zz^\nu\ .\]
    (Here $\Zz^\nu$ denotes the base change to $\MM^2_{GM4}$ of the universal family of EPW sextics over $\MM_{dEPW}$.)
 Possibly after shrinking the base, there is also an isomorphism $\Ss\cong\Ss^\nu$ of the families of surfaces that are the singular locus of each fiber.   
    
The next step is to construct a morphism $\beta$ from the Hilbert scheme of conics $F(Y)$ to the hypersurface $Z\subset I^\vee$. This is done in \cite[Section 4]{IM}; the morphism
$\beta$ is induced (right after Proposition 4.11 in loc. cit.) by a morphism $\alpha$ defined in a geometric way (right after Proposition 4.8 in loc. cit.). Taking the Stein factorization of $\alpha$, one obtains 
    \[ F(Y) \ \to\ X\ \to\ Z\ ,\]
    where $X$ is identified as a double EPW sextic. All of this can be done relative over the base, giving morphisms of families over $\MM^2_{GM4}$
    \[ \FF\ \to\ \XX\ \to \Zz\ . \]
  The argument of \cite[Proposition 4.18]{IM} gives that the isomorphism $\Zz\setminus\Ss\cong \Zz^\nu\setminus\Ss^\nu$ lifts to an isomorphism $\XX\setminus\wt{\Ss}\cong \XX^\nu\setminus \wt{\Ss}^\nu$, where $\wt{\Ss}$ and $\wt{\Ss}^\nu$ denote the fixed loci of the involutions on $\XX$ resp. $\XX^\nu$.
 That is, there is a birational map
   \[ f \colon\ \ \XX\ \dashrightarrow\ \XX^\nu\ \]
   which is, thanks to the results of \cite{IM} (cf. Theorem \ref{im}), a fiberwise isomorphism.
   
   Taking a resolution of indeterminacy of $f$ (resp. of the inverse of $f$), and observing that restriction to a fiber commutes with proper pushforward and with pullback along lci morphisms, this shows that
     \[   f_\ast\bigl( GDA^\ast_{\MM^2_{GM4}}(X)\bigr)\ \subset\ GDA^\ast_{\MM^2_{GM4}}(X^\nu)\ ,\ \ \   f^\ast\bigl( GDA^\ast_{\MM^2_{GM4}}(X^\nu)\bigr)\ \subset\ GDA^\ast_{\MM^2_{GM4}}(X)    \ ,\]
   which closes the proof.
   \end{proof}
  
%
%

  \subsection{A quadratic relation}
  
  \begin{definition}[Zhang \cite{Zh}]  Given a general Gushel--Mukai fourfold $Y$, let $F=F(Y)$ be the variety of conics. We define 
      \[I_F:= \bigl\{ (c,c^\prime)\in F\times F\ \big\vert \ c\cap c^\prime\not=\emptyset\bigr\}\ \ \  \subset F\times F\] 
      as the subvariety given by pairs of intersecting conics.
      
  We also define 
      \[  W:= \bigl\{  (c,c^\prime)   \ \big\vert \exists V_4 \hbox{\ such\ that\  $c$\ and\ $c^\prime$\ are\ residual\ in\ $S_{V_4}$}  \bigr\}\ \ \ \subset F\times F\ ,\]
 where $V_4\subset V_5$ is a 4-dimensional subvector space, and $S_{V_4}$ is the surface $\Gr(2,V_4)\cap H\cap Q$ contained in $Y:=\Gr(2,V_5)\cap H\cap Q$.
   \end{definition}
   
   \begin{remark} Let $P\subset F\times Y$ be the universal conic. The variety $I_F\subset F\times F$ has  dimension 8. As shown in \cite[Lemma 4.7]{Zh}, there is equality
      \begin{equation}\label{Ipq}  I_F = (p\times p)_\ast (q\times q)^\ast (\Delta_Y) = {}^t P\circ P\ \ \ \hbox{in}\ A^2(F\times F)\ .\end{equation}
      The variety $W$ is 6-dimensional.
              \end{remark}
              
   \begin{lemma}[\cite{Lbv}]\label{w} Let $Y$ be a general Gushel--Mukai fourfold, $F=F(Y)$ the variety of conics and $\pi\colon F\to X$ the morphism to the associated double EPW sextic.
   Let $h_F\in A^1(F)$ be an ample divisor, and let $\iota\in\aut(X)$ denote the covering involution.
   There is an equality
     \[ (\pi\times\pi)_\ast \bigl(W\cdot (h_F\times h_F)\bigr)   = d\, \Gamma_\iota + R\ \ \ \hbox{in}\ A^4(X\times X)\ ,\]
     where $d\in\QQ^\ast$ is some constant depending on $h_F$, and $R$ is supported on $x_\rho\times X\cup x_\sigma\times X\cup X\times x_\rho\cup X\times x_\sigma$.
        \end{lemma}  
   
   \begin{proof} This is Lemma 2.23 from \cite{Lbv} and basically follows from the explicit description of $\iota$ in terms of conics given in \cite[Section 3.4]{Zh}.
    \end{proof}

  \begin{proposition}\label{zh}(Zhang \cite{Zh}) Let $Y$ be a general Gushel--Mukai fourfold, and let $F=F(Y)$ be the variety of conics contained in $Y$. There is a relation
    \[  (I_F)^2= \alpha W + I_F\cdot A + B + C\ \ \ \hbox{in}\ A^4(F\times F)\ ,\]
    where $\alpha\in\QQ^\ast$, $A$ and $B$ are decomposable and generically defined (with respect to the family $\XX\to\MM^2_{GM4}$), i.e.
      \[ A,B\ \ \in\ \Bigl\langle   (p_i)^\ast  GDA^\ast_{\MM^2_{GM4}}(F)\Bigr\rangle\ \] 
      (with $p_i$ the projection from $F\times F$ to the $i$-th factor), and $C$ is supported on 
        \[ \Sigma_2:=(F_\rho(Y)\cup F_\sigma(Y))\times (F_\rho(Y)\cup F_\sigma(Y))\ .\]
    \end{proposition}
    
    \begin{proof} See \cite[Proposition 4.10]{Zh} or \cite[Proposition 2.24]{Lbv}. 
     \end{proof}
    
    For our purposes, we will prefer to obtain a quadratic relation on $X\times X$ rather than on $F\times F$ (the reason being that the variety of conics $F$ behaves well only for the {\em general\/} double EPW sextic; another reason being that we want to involve the covering involution $\iota$, which exists on $X$ but a priori not on $F$). To this end, we propose the following definition:
    
   \begin{definition}[\cite{Lbv}] Let $Y$ be a Gushel--Mukai fourfold that is general in the sense of Theorem \ref{im}, let $F=F(Y)$ be the variety of conics, and let $\pi\colon F\to X$ be the morphism to the associated
   double EPW sextic. We define the subvariety $I_\pi\subset X\times X$ by
     \[  I_\pi:= (\pi\times\pi)(I_F)\ \ \ \subset\ X\times X\ .\]
      \end{definition} 
      
      \begin{remark} As observed in \cite[Remark 2.26]{Lbv}, almost all fibers of $\pi\times\pi\colon I_F\to I_\pi$ are of the form $\PP^1\times\PP^1$, and so $I_\pi\subset X\times X$ is 6-dimensional.
      
    We also remark that a priori, the subvariety $I_\pi$ is {\em not \/} intrinsic to $X$. That is, if $Y_1, Y_2$ are two Gushel--Mukai fourfolds with the same associated double EPW sextic $X$
      ($Y_1$ and $Y_2$ are called ``period partners''), it is not clear whether $I_{\pi_1}$ and $I_{\pi_2}$ coincide.
      
%
%
%
        \end{remark}

    \begin{proposition}[\cite{Lbv}]\label{zh2} Let 
    $Y$ be a Gushel--Mukai fourfold that is general in the sense of Theorem \ref{im}, let $X$ be the associated double EPW sextic, and let $\iota\in \aut(X)$ denote the covering involution. There is a relation
      \[ (I_\pi)^2 = \alpha \Gamma_\iota + I_\pi\cdot A + B\ \ \ \hbox{in}\ A^4(X\times X)\ ,\]
      where $\alpha\in\QQ^\ast$, and $A$ and $B$ 
      are as in Proposition \ref{zh}.
     \end{proposition}
   
   \begin{proof}  
     This is obtained in \cite[Proposition 2.27]{Lbv} as a corollary of Zhang's quadratic relation (Proposition \ref{zh}).
      \end{proof}

  \subsection{A linear relation}\label{linearrel}
 
Let us recall the following proposition from \cite[Proposition 2.28]{Lbv}.

  \begin{proposition}\label{I} Let $X$ be a double EPW sextic that is general in the sense of Theorem \ref{im}, and let $\iota$ denote the covering involution. There is a relation
    \[ I_\pi + (\iota,\ide)^\ast (I_\pi) = Z\ \ \ \hbox{in}\ A^2(X\times X)\ ,\]
    where $Z$ is decomposable and generically defined, that is 
      \[  Z\ \ \in\ \Bigl\langle (p_1)^\ast GDA_{\MM^2_{GM4}}^\ast(X), (p_2)^\ast GDA_{\MM^2_{GM4}}^\ast(X)\Bigr\rangle\ .\]
   \end{proposition}

Let us introduce the two correspondences:

\begin{align*}
&\Delta_X^-:= \frac{1}{2}(\Delta_X - \Gamma_\imath)\ ,\\
&\Delta_X^+:= \frac{1}{2}(\Delta_X + \Gamma_\imath)\ \ \ \in A^4(X\times X)\ .
\end{align*}

Hence $A^\ast(X)^-=(\Delta_X^-)_* A^\ast(X)$ and $A^\ast(X)^+=(\Delta_X^+)_* A^\ast(X)$. Before proceeding, let us recall from \cite{Lbv} the following result:

\begin{proposition}\label{deco} Let $X$ be a very general double EPW sextic, and let $\iota$ denote its covering involution. There is a relation
 \[ \Delta_X^- = I \cdot A +B\ \ \ \hbox{in}\ B^4(X\times X)\ ,\]
where $I\in B^2(X\times X)$ is generically defined with respect to the family $\XX\to\MM^2_{GM4}$, and
$A, B$ are in $ \langle (p_i)^*(h), (p_i)^*c_2 \rangle$ (here $p_1, p_2$ denote the projections to the first resp. second factor).
\end{proposition}

\begin{proof} This is essentially contained in the proof of  \cite[Proposition 3.2]{Lbv}. However, as the statement is slightly different from \cite{Lbv}, let us indicate the proof.
By developing the relation given by Proposition \ref{I} 
and taking squares, we get

    \begin{equation}\label{square}  \begin{split} (\iota,\ide)^\ast( I_\pi^2)   &= \bigl( (\iota,\ide)^\ast I_\pi\bigr)^2   \\
                                                                                                &= I_\pi^2-2I_\pi\cdot Z +Z^2\ \ \ \ \hbox{in}\ B^4(X\times X)\ .\\
                                                                                                \end{split} \end{equation}

    Plugging in the relation of Proposition \ref{zh2} into the equality \eqref{square}, we obtain

       \[      (\iota,\ide)^\ast \Bigl(     \alpha \Gamma_\iota + I_\pi\cdot A + B\Bigr) =\alpha \Gamma_\iota + I_\pi\cdot A + B -2I_\pi\cdot Z +Z^2\ \ \ \ \hbox{in}\ B^4(X\times X)\ .\]

Since $(\iota,\ide)^\ast (\Gamma_\iota)$ is the diagonal $\Delta_X$, one sees that this is equivalent to

         \begin{equation}\label{rewrite} \alpha (\Delta_X-\Gamma_\iota) = -(\iota,\ide)^\ast (I_\pi\cdot A) + I_\pi\cdot A +B- 2I_\pi\cdot Z +Z^2 \ \ \ \hbox{in}\ B^4(X\times X)\ .\end{equation}       

  We observe that $A$ and $B$ are decomposable of codimension 2 resp. 4, and that all cycles are generically defined with respect to the family $\XX\to \MM_{GM4}^2$ as in Corollary \ref{Fr2}. 
   In particular, the cycle $Z^2$ is decomposable and generically defined since $Z$ is. What's more, we know (Corollary \ref{Fr2}) that the correspondence $A$ is built from $h\in B^1(X)$ and $c_2\in B^2(X)$ and so $A$ is $\iota$-invariant; hence we have

       \[  (\iota,\ide)^\ast (I_\pi\cdot A) =(\iota,\ide)^\ast (I_\pi)\cdot A\ \ \ \hbox{in}\ B^4(X\times X)\ .\]

     Thus, up to changing $A$ and $B$ (but keeping them generically defined and decomposable), and replacing $I_\pi$ by another generically defined correspondence $I$, we can rewrite \eqref{rewrite}
     as follows
     \[ \alpha (\Delta_X-\Gamma_\iota) =  I\cdot A +B \ \ \ \hbox{in}\ B^4(X\times X)\ , \]
     where $A,B$ are decomposable and generically defined (with respect to the family $\XX\to \MM_{GM4}^2$). Corollary \ref{Fr2} then guarantees that $A,B$ are of
      the required form.
     \end{proof}

\subsection{Voisin's filtration $S_\ast A_0(M)$}
    
Let $M$ be a hyper-Kähler variety de dimension $n=2m$, for zero cycles Voisin \cite{Vcoiso} has constructed a filtration $S_\ast A_0(M) \subset A_0(M)$ in \cite{Vcoiso}, defined as

$$S_i A_0(M):=  \langle x\in M| \dim O_x\geq m-i  \rangle\ ,$$
where $O_x$ is the orbit of $x$ under rational equivalence, and $\langle x \rangle$ indicates the subgroup generated by the points $x$. This is not exactly the definition given in \cite{Vcoiso}; we use the opposite convention, as has also been done for example in \cite{SYZ}, \cite{LZ}, \cite{vialbirat}, \cite{Zh}.

We need to recall a result from \cite[Theorem 1.6]{CMP}, with lighter hypotheses, according to \cite[Theorem 1.5]{MP}:

\begin{theorem}[Charles--Mongardi--Pacienza \cite{CMP}]\label{chmopa}
Let $X$ be a hyper-K\"ahler  variety of $K3^{[m]}$ or generalized Kummer type, and $\dim X=2m$. Then for every non-torsion $L\in \pic (X)$, we have

\begin{align}
S_{m-1}A_{hom}^{2m}(X)_{} & =  L\cdot A_{hom}^{2m-1}(X)_{}\ ; \label{uno} \\
  S_{m-1}A^{2m}(X) & =  L\cdot A^{2m-1}(X)\ . \label{due}
\end{align}
\end{theorem}

For later use, we observe the following:

\begin{corollary}\label{birat} Let $\phi\colon X \dashrightarrow X^\prime$ be a birational map, where $X, X^\prime$ are hyper-K\"ahler varieties of  $K3^{[m]}$-type or generalized Kummer type, and $\dim X=2m$. The isomorphism $\phi_\ast\colon A^{2m}(X)\cong A^{2m}(X^\prime)$ induces an isomorphism
  \[  \phi_\ast\colon\ \ S_{m-1} A^{2m}(X)\ \xrightarrow{\cong}\ S_{m-1} A^{2m}(X^\prime)\ .\]
\end{corollary}

\begin{proof} Thanks to Theorem \ref{chmopa}, it will suffice to show the inclusion
  \begin{equation}\label{preserve}  \phi_\ast \bigl(  A^1(X)\cdot A^{2m-1}(X)\bigr) \ \subset\ \bigl(  A^1(X^\prime)\cdot A^{2m-1}(X^\prime)\bigr)\ .\end{equation}
  
  We will give two proofs for this inclusion. The first proof uses the work of Rie{\ss} \cite[Theorem 3.2]{Rie}, where it is proven that there exists a correspondence $R_\phi$ inducing an isomorphism of graded rings
    \[ (R_\phi)_\ast\colon\ \ A^\ast(X)\ \xrightarrow{\cong}\ A^\ast(X^\prime)\ .\]
    Since it is known that $(R_\phi)_\ast$ coincides with $\phi_\ast$ on 0-cycles \cite[Theorem 3.3(\rom3)]{LV}, this proves the inclusion \eqref{preserve}.
  
For the second (and more direct) proof, let $Z$ denote a resolution of indeterminacy, with morphisms $p\colon Z\to X$ and $q\colon Z\to X^\prime$ such that 
  \[ \phi_\ast=q_\ast p^\ast\colon A^\ast(X)\to A^\ast(X^\prime)\ .\]
Let $E\subset Z$ be the exceptional locus. Since $\phi$ is an isomorphism in codimension 1, we have short exact sequences
  \[ \begin{split}  0 &\to A^{2m-1}(X) \to A^{2m-1}(Z) \to A^{2m-1}(E)\ \\
                          0 &\to A^{2m-1}(X^\prime) \to A^{2m-1}(Z) \to A^{2m-1}(E)\ \\    
                          \end{split}\]
               (Here, $A^\ast(E)$ denotes operational Chow cohomology in the sense of Fulton \cite{F}, and the exact sequences have been constructed by Kimura \cite{Ki}, cf. also \cite{GS}). It follows that there is equality
               \[  p^\ast A^{2m-1}(X) = q^\ast A^{2m-1}(X^\prime)\ \ \ \hbox{in}\ A^{2m-1}(Z)\ .\]
Now, take a curve class $C\in A^{2m-1}(X)$ and a divisor $D\in A^1(X)$. There exists $C^\prime\in A^{2m-1}(X^\prime)$ such that $p^\ast(C)=q^\ast(C^\prime)$, and so
  \[            \begin{split} \phi_\ast (C\cdot D) &= q_\ast p^\ast (C\cdot D)\\
                                         &= q_\ast \bigl(p^\ast(C)\cdot p^\ast(D)\bigr)  \\
                                         &= q_\ast \bigl( q^\ast(C^\prime)\cdot p^\ast(D)\bigr)\\
                                         &= C^\prime\cdot q_\ast p^\ast(D)\ \ \ \hbox{in}\ A^{2m}(X^\prime)\ ,\\
                                         \end{split}\]
proving inclusion \eqref{preserve}.                                                      
  \end{proof}
  
  \begin{remark} It is to be expected that a birational map between hyper-K\"ahler varieties preserves $S_i A^{2m}()$ for all $i$. We don't know how to prove this for $m-i>1$.
   \end{remark}

\section{Main result}

The goal of this section is to prove the following theorem, which is our main result:

\begin{theorem}\label{mainequality}
Let $p:X\to Z$ be a double EPW sextic, and $\imath$ the anti-symplectic involution on $X$. Then we have the equality

$$ A^4(X)^-= S_1 A^4(X) \cap A^4_{hom}(X)\ .$$
\end{theorem}

%
%

\begin{proof}
By equality \eqref{due} of Theorem \ref{chmopa},  we have 
 \[ S_1A^4(X)=A^1(X)\cdot A^3(X)\ ,\] 
 where by this we mean the zero-cycles obtained as sums of intersections of a 1-cycle with a 3-cycle in the Chow ring. Moreover, it is easily seen that
   \[ A^4(X)^-\subset A^4_{hom}(X)_{}\ \]
  (indeed, for any $w\in A^4(X)^-$, we have $p_*(w)=0$ in $A^4(Z)$).   
  Hence, in order to show the inclusion $\subseteq$ it will be enough to show
 that $A^4(X)^- \subset A^1(X)\cdot A^3(X)$.
 To this end, the main ingredient is the following relation of correspondences:
 
 \begin{proposition}[\cite{Lbv}]\label{decodeltaminus}
Let $X$ be any smooth double EPW sextic, and let $\iota$ denote its covering involution. There is a relation

\begin{equation}\label{decomposition}
\Delta_X^-= B + \mathcal{J}_1 + \dots + \mathcal{J}_r \ \ \hbox{in}\ A^4(X\times X)\ ,
\end{equation} 
with the following properties.

\begin{itemize}
\item $B$ is decomposable and generically defined, i.e.
  $$ B\in \bigl\langle  (p_1)^\ast GDA^\ast_{\MM_{dEPW}}(X),  (p_2)^\ast GDA^\ast_{\MM_{dEPW}}(X)\bigr\rangle\ ;$$
\item each $\mathcal{J}_i$ is a composition of the correspondences $\Delta_X^-$ and $I\cdot A$, where

\begin{align*}
I & \in  GDA^2_{\MM_{dEPW}}(X\times X)\ ,\\
A & \in   \bigl\langle (p_1)^*(h), (p_2)^\ast(h), (p_1)^\ast(c_2), (p_2)^\ast(c_2)\bigr\rangle\ , \\
\end{align*}
with $I\cdot A$ occurring at least once in each ${\mathcal J}_i$.

\end{itemize}    
    \end{proposition}
    
\begin{proof}
This is essentially \cite[Proposition 3.2]{Lbv}. However, as the statement is somewhat different from \cite{Lbv}, we indicate the proof.

Since all cycles in the equality are generically defined (with respect to $\MM_{dEPW}$), thanks to the spread lemma (Lemma \ref{spread}), it is enough to prove the
equality of Proposition \ref{decodeltaminus} for a very general double EPW sextic $X$.
Thus, let us now consider the universal family of double EPW sextics $\XX\to \MM_{dEPW}$, and a very general fiber $X$. 

    
   Recall that there is a morphism $\nu\colon \MM^2_{GM4}\to \MM_{dEPW}$ (Theorem \ref{im}).                              
          Cutting with generic hyperplane sections $h_i$ and shrinking $\MM^2_{GM4}$, we obtain a finite degree $d$ morphism
          \[ \nu\colon \MM^3_{GM4}:= \MM^2_{GM4}\cap h_1  \cap \cdots\cap h^4\ \to\ \MM_{dEPW}\ .\]
          This gives rise to a fibre diagram
        \[ \begin{array}[c]{ccc}   \XX^\nu & \xrightarrow{\rho} &   \XX \\
                                                     &&\\
                                                     \downarrow && \downarrow\\
                                                     \MM^3_{GM4} & \xrightarrow{\nu} & \MM_{dEPW}\ , \\
                                                     \end{array} \]
   where we define $\XX^\nu:= \XX\times_{\MM_{dEPW}} \MM^3_{GM4}$ as the base changed family.
                                  
   Over the (very general) fiber $X$ of the family $\XX\to \MM_{dEPW}$, there are $d$ (isomorphic) fibers $X_1,\ldots,X_d$ of the family $\XX^\nu\to \MM^3_{GM4}$ (these $X_i$
   correspond to different GM fourfolds with the same associated double EPW sextic). 
   For each of the $X_i$, Proposition \ref{deco} gives a relation of correspondences modulo algebraic equivalence                              
   \begin{equation}\label{eachi}  \Delta_{X_i}^- = I_i \cdot A +B\ \ \ \hbox{in}\ B^4(X_i\times X_i)\ \ \ (i=1,\ldots, d)\ ,   \end{equation}
   where $A, B$ are constant (since they are built out of $h$ and $c_2$) and each $I_i\in B^2(X_i\times X_i)$ is generically defined with respect to the family $\XX\times_{\MM^2_{GM4}} \XX$ of Theorem \ref{im} (more precisely, the proof of Proposition \ref{deco} shows there exists one cycle in $B^2(\XX\times_{\MM^2_{GM4}} \XX)$ restricting to the various $I_i$). 
    Thanks to Lemma \ref{2fam}, $I_i$ is also generically defined with respect to the family $\XX^\nu\times_{\MM^2_{GM4}} \XX^\nu$, and a fortiori with respect to the family $\XX^\nu\times_{\MM^3_{GM4}} \XX^\nu$.
   
   Let us now define a correspondence
     \[ I:= {1\over d}\  \displaystyle\sum_{i=1}^d  \bigl((\rho,\rho)\vert_{X_i\times X_i}\bigr){}_\ast (I_i)\ \ \in\ B^2(X\times X)\ .\]
  Equality \eqref{eachi} ensures that this correspondence $I$ verifies the equality
    \begin{equation}\label{onlyalg}  \Delta_{X}^- = I \cdot A +B\ \ \ \hbox{in}\ B^4(X\times X)\   .\end{equation}
   Moreover, the above considerations guarantee that $I$ (and hence each term in equality \eqref{onlyalg}) is generically defined with respect to $\MM_{dEPW}$.

The final step consists in applying the celebrated nilpotence theorem of Voevodsky \cite{Voe} and Voisin \cite{V2}. This transforms the algebraic equivalence \eqref{onlyalg} into an equality modulo rational equivalence
   \[  \Bigl(   \Delta_X^- - I\cdot A - B    \Bigr)^{\circ M}=0\ \ \ \hbox{in}\ A^4(X\times X) \]
 for some $M\in\NN$. Developing this expression
    (and observing that $ \Delta_X^-$ is idempotent while any generically defined correspondence composed with $B$ is decomposable and generically defined), 
  one obtains the equality of Proposition \ref{decodeltaminus} for the very general fiber $X$ of the family $\XX\to \MM_{dEPW}$.
  \end{proof}

Armed with Proposition \ref{decodeltaminus}, we now continue and prove the inclusion ``$\subset$'' of the theorem. 

Take $a\in A^4(X)^-$. To show that $a\in S_1 A^4(X)= A^1(X)\cdot A^3(X)$, we let both sides of equality \eqref{decomposition} act on $a$. The left-hand side acts as the identity on $a$.
As observed above, $a\in A^4_{hom}(X)$, and hence we have $B_*(a)=0$ since $B$ is a decomposable correspondence. Thus, we need to show that 
  \[ (\mathcal{J}_i)_*(a)\ \ \in\ A^1(X)\cdot A^3(X)\ \ \forall i\ .\] 
We will do this case by case, depending on the different forms the cycle $A$ can take. 
\smallskip

 Hence, we will need to give an argument for the case that $A$ equals $p_1^*(C)$ or $p_2^*(C)$, with $C\in A^2(X)$ a linear combination of $h^2$ and $c_2(X)$, and for the case that $A=p_1^*(h)\cdot p_2^*(h)$, with $h\in A^1(X)$. 

\smallskip

$(i)$ First, let us assume that ${\mathcal J}_i$ contains $I\cdot A$ where
$A$ is of the form $p_1^*(h)\cdot p_2^*(h)$. The projection formula then gives 
    \[(I\cdot A)_*(a) = h\cdot I_*(a\cdot h)\ ,\] 
    and $(a\cdot h)$ is zero for dimension reasons. It follows that $({\mathcal J}_i)_\ast(a)$ is zero.

 $(ii)$ In the second case, let us assume that ${\mathcal J}_i$ contains $I\cdot A$ where
$A$ is of the form $A=p_1^*(C)$. Then, the projection formula gives 
\[(I\cdot A)_*(a)=I_*(a\cdot C)\ ,\] 
which is zero for dimension reasons. Again, it follows that $({\mathcal J}_i)_\ast(a)$ is zero.

$(iii)$ In the third and final case,
let us assume that the correspondence ${\mathcal J}_i$ is a composition of $\Delta_X^-$ and $I\cdot A$, where
  $A=p_2^*(C)$. Left-composing both sides of the equation \eqref{decomposition} with $\Delta_X^-$, we may assume ${\mathcal J}_i$ is of the form
    \[     {\mathcal J}_i=\Delta_X^-\circ (I\cdot (p_2)^\ast(C))\circ {\mathcal J}_i^\prime\ ,\]
    for some ${\mathcal J}_i^\prime\in A^4(X\times X)$. We can thus write 
   \[  ({  \mathcal J}_i)_\ast(a) =  (\Delta_X^-)_\ast (I\cdot (p_2)^\ast(C))_\ast (a^\prime)\ ,\]
   where $a^\prime:= ( {  \mathcal J}^\prime_i)_\ast(a)$.
  Applying the projection formula, this gives 
  \[    ({  \mathcal J}_i)_\ast(a) =  (\Delta_X^-)_\ast \bigl((I_*a^\prime)\cdot C\bigr)\ .\] 
  
  We observe that $C\in A^2(X)^+$, and so this can be rewritten
    \[ \begin{split} ({  \mathcal J}_i)_\ast(a) &=  (\Delta_X^-)_\ast \bigl((I_*a^\prime)\cdot C\bigr)\\
                                                          &= {1\over 2}\Bigl((I_*a^\prime)\cdot C  - \iota_\ast \bigl((I_*a^\prime)\cdot C\bigr)\Bigr)\\
                                                          &=  {1\over 2}\bigl(   (I_*a^\prime)\cdot C  - \iota_\ast (I_*a^\prime)\cdot \iota_\ast(C)\bigr)\\
                                                          &= {1\over 2}\bigl( (I_*a^\prime)\cdot C  - \iota_\ast (I_*a^\prime)\cdot C\bigr)\\     
                                                          &=        \bigl((\Delta_X^-)_\ast (I_*a^\prime)\bigr)\cdot C\\
                                                          &= b\cdot C \  ,\\
                                                          \end{split}\] 
with $b:=  (\Delta_X^-)_\ast (I_*a^\prime) \in A^2(X)_{}^-$.  The cycle $b$ being anti-invariant, it follows that $b$ has zero intersection with the fixed surface $S$ of the involution inducing the double cover $p:X\to E$ (see \cite[Remark 2.7]{LV}). From the work of Ferretti \cite[Lemma 4.1]{Fe}, we know that $S$ is a rational combination (with non-zero coefficients) of $c_2$ and $h^2$ in $A^2(X)$. We thus find that $b\cdot c_2 $ is proportional to $b\cdot h^2$. Since $A=p_2^*(C)$ and, by Proposition \ref{decodeltaminus}, $A\in\langle (p_i)^*(h^2),(p_i)^*c_2 \rangle$  the cycle $C$ is then made of $c_2$ and $h^2$,  
   and so 
   \[ ({\mathcal J}_i)_\ast(a)=b\cdot C =b\cdot h^2\ \ \in A^2(X)^-\cdot h^2 \subset A^1(X)\cdot A^3(X)\ .\]
   
  Combining these 3 cases, we have now proven the inclusion
      \begin{equation}\label{h2} A^4(X)^- \ \subset\   A^2(X)^-\cdot h^2\ \subset\ A^1(X)\cdot A^3(X)\ .\end{equation}

\smallskip

Let us now prove the opposite inclusion $ A^4(X)^-\supset S_1 A^4(X) \cap A^4_{hom}(X)$. Recall that $X$ is a double cover of an EPW sextic $Z\subset\PP^5$. 
It is easy to see that
\begin{equation}\label{1to3}
A^1(Z)\cdot A^3(Z)= \QQ 
\end{equation}
(indeed, $A^1(Z)\cong\QQ$ is spanned by the hyperplane class, and so $A^1(Z)\cdot A^3(Z)=j^\ast j_\ast A^3(Z)=j^\ast A^4(\PP^5)$ where $j\colon Z\hookrightarrow \PP^5$ denotes the inclusion morphism).

Equation (\ref{1to3}) implies that 
  \begin{equation}\label{1to3t} A^1(X)^+\cdot A^3(X)^+ \cap A^4_{hom}(X)=0\ .\end{equation}
  
  In view of equality \eqref{uno} from Theorem \ref{chmopa}, plus the fact that for any double EPW sextic one has $A^1(X)^+=\QQ[h]$ (where $h$ denotes the polarization), we have
  $$S_1A^4_{hom}(X)_{}  =  h\cdot A^3_{hom}(X)_{}\ .$$
It follows that
  \[ \begin{split} S_1A^4_{hom}(X)_{} &=  h \cdot \Bigl( A^3_{hom}(X)^+\oplus A^3_{hom}(X)_{}^-\Bigr)\\
                                  &= h\cdot A^3_{hom}(X)_{}^-\\  
                                   & = A^1(X)^+\cdot A^3_{hom}(X)_{}^-\subset A^4(X)^-\ \\
                                   \end{split}
                                   \] 
                                   (where the second equality is an application of the vanishing \eqref{1to3t}),
  and so the wanted inclusion is proven.

\end{proof}

\begin{corollary}\label{motivicS1} 
The component $S_1A^4(X) = A^1(X)\cdot A^3(X)$ of the Voisin filtration is motivic (that is, it is cut out by a projector).
\end{corollary}     

\begin{remark}\label{stuff}

\noindent
\begin{enumerate}
\item Zhang \cite{Zh} has shown that the surface of fixed points of the involution of $X$ is a constant cycle surface. The class of this surface in $A^2(X)$ is a combination of  $c_2(X)$ and $h^2$ \cite{Fe}. As shown by Voisin \cite[Proof of Lemma 3.10]{Vcoiso}, this implies that $S_0 A^4(X)$ is one-dimensional. Hence,
 $S_0A^4(X)= \QQ [o_X]$ is also motivic.
 
\item Corollary \ref{motivicS1} can be rephrased in other terms by saying that the birational Chow--Künneth decomposition of \cite[Remark 4.8]{vialbirat} respects $S_1A^4(X)$.
An interesting question is whether this birational Chow--K\"unneth decomposition is co-multiplicative, in the sense of \cite{vialbirat}. We did not succeed in settling this.
\item Partial results closely related to our Theorem \ref{mainequality} were obtained in \cite{LV}, see Theorem 3.6 and Corollary 3.12 therein.
\item for a very general double EPW sextic, the inclusion $A^4(X)^-\subseteq S_1A^4(X)\cap A^4_{hom}(X)$ was already proven in \cite[Proposition 3.12]{LYZ} by different methods, building on results of \cite{Zh}. Our proof shows the inclusions in both directions for all double EPW sextics.
\item  The ``deepest piece'' $A^4_{hom}(X)^+$ is still something mysterious. We expect that (similar to known results for the Fano variety of lines on cubic fourfolds \cite[Theorem 4]{SV}, and in accordance with the fact that $H^4(X,{\mathcal O}_X)=\sym^2 H^2(X,{\mathcal O}_X)$) there is equality
  \[ A^4_{hom}(X)^+  = A^2_{hom}(X)^- \cdot A^2_{hom}(X)^-\ . \]
  We have not been able to prove this.
   \end{enumerate}     
\end{remark}







\section{Some consequences}

\subsection{Chow--Lefschetz}

\begin{proposition}\label{lefschetz} Let $X$ be any double EPW sextic.

\noindent
(\rom1) There are isomorphisms
  \[ \begin{split} &\cdot h\colon\ \ A^2_{hom}(X)^-\ \xrightarrow{\cong}\ A^3_{hom}(X)^-\ ,\\
                          & \cdot h\colon\ \ A^3_{hom}(X)^-\ \xrightarrow{\cong}\ A^4(X)^-\ .\\
                      \end{split}       \]
  
  \noindent
  (\rom2) There is a surjection
    \[ \cdot h^2\colon\ \ A^2(X)\ \twoheadrightarrow\ S_1 A^4(X)\ .\]
    \end{proposition}

\begin{proof}
For (\rom1), the surjectivity of the second arrow was already proven in \eqref{h2}. To prove injectivity of the second arrow, let $\gamma\in A^3_{hom}(X)^-$ be such that $\gamma\cdot h=0$ in $A^4(X)$. Applying the equality of correspondences of Proposition \ref{decodeltaminus} to $\gamma$, we find that
  \[ \gamma =  \bigl( {\mathcal J}_1+\cdots +{\mathcal J}_r\bigr){}_\ast(\gamma)   \ \hbox{in}\ A^3(X)\ ,\]  
  for certain correspondences ${\mathcal J}_i$ built from $\Delta_X^-$ and $I\cdot A$. Here $A\in A^2(X\times X)$ is a cycle of the form 
  $A= (p_1)^\ast(a)$ or  $A=(p_2)^\ast(a^\prime)$ or $A= (p_1)^\ast(D)\cdot (p_2)^\ast(D^\prime)$
 where $a,a^\prime\in A^2(X)$ are linear combinations of $h^2$ and $c_2$, and $D, D^\prime\in A^1(X)$.
 
 Let us do a case by case analysis, similar to the proof of Theorem \ref{mainequality}.
In case ${\mathcal J}_i$ is a correspondence containing a factor $I\cdot (p_1)^\ast(a)$ or a factor $I\cdot (p_2)^\ast(a)$, 
 it is easily seen that ${\mathcal J}_i$ acts as zero on $\gamma$ (indeed, the action factors over $A^5(X)=0$ resp. over $A^1_{hom}(X)=0$).
 The only remaining possibility is thus that ${\mathcal J}_i$ is a correspondence which is a composition of $\Delta_X^-$ and $I\cdot  (p_1)^\ast(D)\cdot (p_2)^\ast(D^\prime)$.
 As in the proof of Theorem \ref{mainequality}, up to right-composing both sides of the equality \eqref{decomposition} with $\Delta_X^-$, we may suppose ${\mathcal J}_i$ is of the form
    \[     {\mathcal J}_i=  {\mathcal J}_i^\prime \circ \bigl(I\cdot  (p_1)^\ast(D)\cdot (p_2)^\ast(D^\prime)\bigr)\circ \Delta_X^-\ ,\]
    for some ${\mathcal J}_i^\prime\in A^4(X\times X)$.    
  It follows that we can write 
   \[ ({\mathcal J}_i)_\ast(\gamma)=   ({\mathcal J}_i^\prime)_\ast \bigl( I\cdot (p_1)^\ast(D)\cdot (p_2)^\ast(D^\prime) \bigr){}_\ast(\gamma)\ .\]
 Applying the projection formula, we now find that
   \[  ({\mathcal J}_i)_\ast(\gamma)=   ({\mathcal J}_i^\prime)_\ast \bigl(D^\prime\cdot I_\ast (\gamma\cdot D)\bigr) \ \ \hbox{in}\ A^3(X)\ .\]
   But $D$, being generically defined, is proportional to the polarization $h$ and so we obtain $\gamma\cdot D=0$ and hence
      \[   ({\mathcal J}_i)_\ast(\gamma)=0\ \ \hbox{in}\ A^3(X)\ .\]
This proves injectivity of the second arrow. The bijectivity of the first arrow is proven in a similar fashion.
         
  
To prove (\rom2), in view of (\rom1) and Theorem \ref{mainequality}, it suffices to observe that
  \[ S_1 A^4(X)\cap A^4(X)^+ =\QQ[h^4] = A^2(X)^+\cdot h^2\ .\]
\end{proof}

\begin{remark} One can consider Proposition \ref{lefschetz}(\rom1) as a ``motivic manifestation'' of the hard Lefschetz isomorphisms in cohomology 
  \[ H^2(X,\QQ)^- \xrightarrow{\cong} H^4(X,\QQ)^-\xrightarrow{\cong} H^6(X,\QQ)^-\ .\] 
  
  The equality 
  \[ S_1 A^4(X)= A^2(X)\cdot h^2\]
  of Proposition \ref{lefschetz}(\rom2)
also holds for Fano varieties of lines in cubic fourfolds; this follows from results of Shen--Vial \cite{SV} combined with \cite[Proposition 4.5]{Vcoiso}.
\end{remark}

\subsection{Chow--Abel--Jacobi}

    
\begin{proposition}[Zhang \cite{Zh}]\label{0zh}  
Let $Y$ be a GM fourfold such that there exists an associated double EPW sextic $X$ via the construction of Theorem \ref{im}.  Then we have a correspondence $\Gamma$ inducing an isomorphism

$$  \Gamma_\ast\colon\ \  A^3_{hom}(Y)\ \xrightarrow{\cong}\ A^4(X)^-\ .$$
\end{proposition}
     
\begin{proof} Recall that the universal conic $P$ for a Gushel-Mukai fourfold $Y$ tautologically lives in $F\times Y$, but $F$ has a natural map $\pi$ to the double EPW sextic $X$ associated to $Y$ via the construction of Theorem \ref{im}. Let us denote the two projections by $p:P\to F$ and $q:P\to Y$.
Starting from $ A^3_{hom}(Y)$, we apply the Abel-Jacobi map $p_*q^*$, then intersect with a hyperplane section of $F$, then push forward to $X$. By \cite[Proposition 5.4]{Zh}, this induces an isomorphism onto $A^4(X)^-$.
\end{proof}


%
%

Let us now upgrade Proposition \ref{0zh} to an isomorphism of Chow motives. First, let us recall from \cite[Section 7.4]{FM} that there exists a refined Chow-Künneth decomposition for the Chow motive of a Gushel--Mukai fourfold $Y$, such that $h^4(Y)$ decomposes as $h^4_{alg}(Y)\oplus t^4(Y)$. These are the algebraic and transcendental parts of $h^4$, which realize the algebraic respectively the transcendental parts of the cohomology. For simplicity, from now on we will denote by $t(Y)$ the transcendental motive of a Gushel--Mukai fourfold.

\begin{proposition}\label{isomotives} Let $Y$ be a Gushel--Mukai fourfold such that there exists an associated double EPW sextic $X$ via the construction of Theorem \ref{im}. Let us write $X\to Z$ for the double cover, where $Z\subset\PP^5$ is an EPW sextic. There is an isomorphism of Chow motives
  \[  h(X) \cong h(Z) \oplus \bigoplus_{i=-1}^1 t(Y)(i)\oplus \bigoplus_{j=-3}^{-1} \one(j)^{\oplus r_j}\ \ \ \hbox{in} \ \MM_{\rm rat}\ .\]
\end{proposition}

\begin{proof} It will suffice to establish an isomorphism
  \[ h(X)^-\cong \bigoplus_{i=-1}^1 t(Y)(i)\oplus \bigoplus_{j=-3}^{-1} \one(j)^{\oplus r_j}\ \ \ \hbox{in}\ \MM_{\rm rat}\ ,\]
where $h(X)^-:=(X,\Delta_X^-,0)$ denotes the anti-invariant part of the motive of $X$.
As we work over $\C$ which is a universal domain, the Bloch--Srinivas decomposition argument \cite{BS} (cf. also \cite[Lemma 1.1]{Huy2}) reduces this to showing that there are correspondence-induced isomorphisms of Chow groups
  \[ A^\ast_{}(X)^-\cong  A^\ast_{}\bigl(  \bigoplus_{i=-1}^1 t(Y)(i) \oplus \bigoplus_{j=-3}^{-1}  \one(j)^{\oplus r_j} \bigr)\ .\]
This readily follows from the combination of the Abel--Jacobi isomorphism in cohomology, Zhang's result (Proposition \ref{0zh}), and the isomorphisms of Proposition \ref{lefschetz}
  \[  A^2_{hom}(X)^-\cong A^3_{hom}(X)^-\cong A^4(X)^-\ .\]  
\end{proof}

\begin{remark}  Let $Y_1, Y_2$ be 2 Gushel--Mukai fourfolds that are general in the sense of Theorem \ref{im} and with the same associated double EPW sextic (these are called period partners). It is immediate from Proposition \ref{0zh} or \ref{isomotives} that $Y_1$ and $Y_2$ have isomorphic Chow motives. This was already proven by Fu--Moonen \cite{FM}, but their proof relies on deep results on the level of derived categories obtained in \cite{KP}; the alternative proof given here is direct and geometric.
\end{remark}

\begin{remark} Let $\XX\to \MM_{dEPW}$ denote the universal family. Proposition \ref{isomotives} (plus the Franchetta property for the universal family of Gushel--Mukai fourfolds, which is easy, cf. \cite{BoLa}) readily implies that the anti-invariant part of this family has the Franchetta property, i.e. the cycle class map induces injections
  \[  GDA^\ast_{\MM_{dEPW}}(X)^- \ \hookrightarrow\ H^\ast(X,\QQ)\ .\] 
  \end{remark}

\subsection{Cubic fourfolds in ${\mathcal C}_{12}$}

The goal of this section is to show that Conjecture \ref{conj} holds true for the  Fano varieties of lines $X=F(Y)$ of smooth cubic fourfolds $Y$ contained in the Hassett divisor $\mathcal{C}_{12}$. This divisor is the closure of the locus of cubics containing a smooth cubic scroll. Hassett and Tschinkel have described \cite{HT}
the ample and moving cone in $\pic (X)$ and described birational automorphisms of infinite order in $\bir(X)$.

In \cite{BFMQ}, an interesting theory of birational models of Fano varieties of lines of cubics in $\mathcal{C}_{12}$ is developed. The most interesting aspect with respect to this paper is that the presence of a cubic scroll inside the cubic fourfold influences in an important way the birational geometry of the Fano variety of $Y$, as shown by the following result.

\begin{theorem}[Brooke--Frei--Marquand--Qin \cite{BFMQ}]\label{birmodels}
Let $X=F(Y)$ be the Fano variety of lines of a very general cubic fourfold $Y \in \mathcal{C}_{12}$. Then $X$ has three isomorphism classes of birational hyper-Kähler models, represented by $X$ itself and two non-isomorphic Mukai flops, both of which are isomorphic to double EPW sextics. Moreover, $\bir(X)$ is generated by the covering involutions on these two double EPW sextics.
\end{theorem}

\begin{proof} This is \cite[Theorem 1.1]{BFMQ}.
\end{proof}

\begin{proposition}\label{conjC12}
Let $Y$ be a very general cubic fourfold in the divisor $\mathcal{C}_{12}$, let $X=F(Y)$ be its Fano variety of lines and let $\phi\in\bir(X)$ be any birational automorphism. Then Conjecture \ref{conj} holds true for $(X,\phi)$. 
\end{proposition}

\begin{proof}
We know that $S_0 A^4(X)=\QQ[c_4(X)]$ (this is true for the Fano variety of lines in any cubic fourfold, cf. \cite[Proposition 4.5]{Vcoiso}), and so any birational automorphism acts as the identity on $S_0 A^4(X)$. It remains to prove that if $g\in\bir(X)$ is one of the two covering involutions that generate $\bir(X)$ (hence $g$ is anti-symplectic) and $A^4(X)^-$ denotes the $(-1)$-eigenspace with respect to $g$ then
  \[ A^4(X)^- = S_1 A^4(X)\cap A^4_{hom}(X)\ .\]

By Theorem \ref{birmodels}, we know that $g$ comes from the covering involution of some double EPW sextic $W$ birational to $X$. Now we observe that, by Corollary \ref{birat}, we can consider the action of $\bir(X)$ equivalently on $S_1A^4(X)$ or on $S_1A^4(W)$. But, thanks to our Theorem \ref{mainequality}, we know the above equality holds true for $S_1A^4(W)$ and its covering involution.
\end{proof} 

\begin{proposition} Let $Y$ and $X$ be as in Proposition \ref{conjC12}, and let $\phi\in\bir(X)$ be any anti-symplectic birational automorphism. Let $A^\ast(X)^-$ denote the $(-1)$-eigenspace with respect to the action of $\phi$ on Chow groups. There are correspondence-induced isomorphisms
  \[  A^2_{hom}(X)^-\cong A^3_{hom}(X)^-\cong A^4(X)^-\cong A^3_{hom}(Y)\ .\]
  \end{proposition}
  
 \begin{proof} Without loss of generality, we may assume $\phi$ is a generator of $\bir(X)$.
 It follows from Theorem \ref{birmodels} that there exists a double EPW sextic $W$ birational to $X$ such that $\phi$ is induced by the covering involution of $W$.
 There are isomorphisms of Chow groups $A^i(X)\cong A^i(W)$ sending $A^i_{hom}(X)^-$ to $A^i_{hom}(W)^-$, and vice versa
 (in codimension $i=3,4$, this is because the ring isomorphism of \cite{Rie} coincides with the map given by the closure of the graph of the birational map, cf. \cite[Theorem 3.3(\rom3)]{LV}; in codimension $i=2$, it can be seen directly that the closure of the graph of the birational map induces an isomorphism $A^2_{hom}(X)\cong A^2_{hom}(W)$).
 
This implies there are isomorphisms
  \[ A^i_{hom}(X)^-\cong A^i_{hom}(W)^-\ \ \ (i=2,3,4)\ ,\]
  where the right-hand side refers to the covering involution of $W$. The Chow--Lefschetz isomorphisms of Proposition \ref{lefschetz} then give isomorphisms
    \[ A^2_{hom}(X)^-\cong A^3_{hom}(X)^-\cong A^4(X)^-\ .\]
    On the other hand, we also have
    \[ A^4(X)^-= S_1 A^4(X)\cap A^4_{hom}(X) = A^4_{(2)}(X) \ ,\]
    where $A^\ast_{(\ast)}(X)$ refers to the Shen--Vial splitting \cite{SV}; here the first equality is Proposition \ref{conjC12} while the second equality is part of \cite[Proposition 4.5]{Vcoiso}. Since $A^4_{(2)}(X)$ is isomorphic to $A^3_{hom}(Y)$ via the universal line \cite[Proposition 19.5 and Theorem 20.5]{SV}, this closes the proof.
  \end{proof}

\begin{remark}
In \cite[Theorem 1.2 and Theorem 1.3]{BFMQ}, two other (18-dimensional) families of cubics $Y$ in $\mathcal{C}_{12}$ with $rk(H^{2,2}(Y)\cap H^4(Y,\mathbb{Z}))> 2$ are introduced. These families have respectively five and eight isomorphism classes of birational hyper-Kähler models, and most of these models are double EPW sextics. For the very general element in the family, $\bir(X)$ is again generated by the covering involutions of the double EPW sextics \cite[Theorem 1.5]{BFMQ}.
Hence the argument of Proposition \ref{conjC12} shows that Conjecture \ref{conj} holds true for these families as well.
\end{remark}

\vskip0.5cm

 \begin{nonumberingt} RL thanks MB for an invitation to sunny Montpellier (June 2024), where this project originated.
MB and RL are supported by the ANR project FanoHK (ANR-20-CE40-0023). MB is a member of the Réseau thématique \rm Géométrie algébrique et singularités, \it and of \rm  GNSAGA. \it Finally, we are very grateful to the referee for lots of constructive comments and corrections, saving us from embarrassment in several places.
\end{nonumberingt}

\end{document}